\documentclass[11pt]{amsart}
\usepackage{amsmath,amsthm, amscd, amssymb, amsfonts}
\usepackage[dvips, dvipsnames, usenames]{color}
\usepackage{hhline}

\newcommand\Ur{\mathcal{UR}}

\newcommand\boxe{\begin{tabular}{|p{0,1cm}|}
\hline \\ \hline \end{tabular}}

\newcommand\boxu{\begin{tabular}{|p{0,1cm}|}
\\ \hline  \end{tabular}}

\newcommand\boxa{\begin{tabular}{p{0,1cm}|}
\hline \\ \hline \end{tabular}}

\newcommand\boxur{\begin{tabular}{p{0,1cm}|}
\hline \\  \end{tabular}}

\DeclareMathOperator*{\Tim}{\times}
\newcommand{\Times}[2]{\sideset{_#1}{_#2}\Tim}

\newcommand{\proy}{\varPi}

\newcommand{\acts}{{\rightharpoondown}}

\newcommand{\K}{{\mathcal K}}

\newcommand{\G}{{\mathcal G}}

\newcommand{\Q}{{\mathcal Q}}
\newcommand{\F}{{\mathcal F}}

\newcommand{\D}{{\mathcal D}}
\newcommand{\Ec}{{\mathbf E}}
\newcommand{\Ee}{{\mathcal E}}
\newcommand{\Kc}{{\mathbf K}}

\newcommand\idh{\id_{\Hc}}
\newcommand\idv{\id_{\Vc}}

\newcommand{\B}{{\mathcal B}}
\newcommand{\T}{{\mathcal B}}
\newcommand{\Hc}{{\mathcal H}}
\newcommand{\Vc}{{\mathcal V}}
\newcommand{\Pc}{{\mathcal P}}
\newcommand{\Oc}{{\mathcal O}}

\newcommand{\GQ}{{\mathcal G}_{\curvearrowright Q}}
\newcommand{\Gpx}{{\mathcal G}_{\rightarrow p(x)}}

\newcommand{\path}{\operatorname{Path}}

\newcommand{\fde}{{\triangleright}}

\newcommand{\prin}{t}
\newcommand{\fin}{b}
\newcommand{\pri}{r}
\newcommand{\fine}{l}

\numberwithin{equation}{section}\theoremstyle{plain}

\newtheorem{theorem}{Theorem}[section]
\newtheorem{mtheorem}{Theorem}
\newtheorem{lema}[theorem]{Lemma}

\newtheorem{prop}[theorem]{Proposition}

\newtheorem{zulema}[theorem]{Sublemma}
\newtheorem{zulma}[theorem]{Restrictions}

\theoremstyle{definition}
\newtheorem{definition}[theorem]{Definition}
\newtheorem{exa}[theorem]{Example}

\theoremstyle{remark}
\newtheorem{obs}[theorem]{Remark}

\newtheorem{fc}[theorem]{Filling condition}

\newcommand\id{\operatorname{id}}
\newcommand\idd{\mathbf{id}}
\newcommand\iddh{\mathbf{id}}
\newcommand\iddv{\mathbf{id}}
\newcommand\op{\operatorname{op}}

\def\pf{\begin{proof}}
\def\epf{\end{proof}}

\theoremstyle{remark}

\begin{document}

\renewcommand{\baselinestretch}{1.2}

\thispagestyle{empty}
\title[Double groupoids]{The structure of double groupoids}
\author[Andruskiewitsch and Natale]{ Nicol\'as Andruskiewitsch and Sonia Natale}
\address{\noindent
Facultad de Matem\'atica, Astronom\'\i a y F\'\i sica,
Universidad Nacional de C\'ordoba.  CIEM -- CONICET.
(5000) Ciudad Universitaria, C\'ordoba, Argentina}
\email{andrus@mate.uncor.edu, \newline \indent \emph{URL:}\/
http://www.mate.uncor.edu/andrus} \email{natale@mate.uncor.edu,
\newline \indent \emph{URL:}\/ http://www.mate.uncor.edu/natale}
\thanks{This work was partially supported by CONICET, Fundaci\' on
Antorchas, Agencia C\'ordoba Ciencia, ANPCyT    and Secyt (UNC)}
\date{October 18, 2007}
\dedicatory{A la memoria de Gustavo Isi, de La Floresta
(1961-2006)}
\begin{abstract} We give a general description of the structure of a
discrete double groupoid  (with an extra, quite natural, filling
condition) in terms of groupoid factorizations and groupoid
$2$-cocycles with coefficients in abelian group bundles. Our
description goes as follows: To any double groupoid, we associate
an abelian group bundle and a second double groupoid, its frame.
The frame satisfies that every box is determined by its edges, and
thus is called a `slim' double groupoid. In a first step, we prove
that every double groupoid is obtained as an extension of its
associated abelian group bundle by its frame.  In a second,
independent, step we prove that every slim double groupoid with
filling condition is completely determined by a factorization of a
certain canonically defined `diagonal' groupoid.
\end{abstract}

\maketitle

\section*{Introduction}

The main result of this paper is the determination of the
structure of a discrete double groupoid -satisfying a natural
filling condition- in terms of groupoid data. By 'discrete' we
mean here that no additional structure (differential, measurable,
topological, etc.) is assumed. The problem of describing all
double groupoids in terms of more familiar structures was explicitly raised
by Brown and Mackenzie in \cite[p. 271]{BM}.

Double groupoids  were introduced by Ehresmann \cite{ehr} in the
early sixties, and later studied by several people because of
their connection with different areas of mathematics, such as
homotopy theory,  differential geometry and Poisson-Lie groups.
See for instance \cite{brown, bj, BM, bs, l, wl, mk1,  mk3, mk2,
pr, pradines} and references therein.

A double groupoid is a groupoid object in the category of
groupoids. This can be interpreted as a set of 'boxes' with two
groupoid compositions --the \emph{vertical} and \emph{horizontal}
compositions--, together with compatible groupoid compositions of
the edges, obeying several conditions, in particular  and most
importantly the so called \emph{interchange law}.

The double groupoids we are interested in satisfy the following:
\begin{fc}\label{filling}For every configuration of matching edges
as in Figure 1 (a), there is at least one box in the double
groupoid as depicted in Figure 1 (b). This box is called the
'filling' of the given configuration.\end{fc}
\begin{equation*}\underset{\text{Figure 1 (a)}}{\begin{matrix}\quad x \quad
\\ \boxur \, g \\ \quad \end{matrix}} \qquad \qquad \qquad \underset{\text{Figure 1 (b)}}{\begin{matrix} \quad x \quad \\ \, \,\,
\boxe \,\, g\\ \quad
\end{matrix}}\end{equation*}
This  condition is often assumed in the case of double groupoids
arising in differential geometry, and is discussed by Mackenzie in
\cite{mk2}.

Concerning the structure of double groupoids, some very early
results on `special double groupoids with special connections'
were obtained in \cite{bs}. For more general double groupoids,
only those in a few classes were known to be determined by
`groupoid data'. One of them is the class of \emph{vacant} double
groupoids ({\it i.e.} those for which every configuration as in
Figure 1 (a) has a unique filling): it was proved by Mackenzie in
\cite{mk1, mk2}, that vacant double groupoids are essentially the
same thing as \emph{exact} factorizations of groupoids. Another
one is the class of \emph{transitive} or \emph{locally trivial}
double groupoids: in the paper \cite{BM}, R. Brown and Mackenzie
show that such a double groupoid is determined by its \emph{core
diagram}.

\medbreak We now state our main result; unexplained notions and
notations will be properly defined later in the text. Let
$\begin{matrix} \B &\rightrightarrows &\Hc
\\\downdownarrows &&\downdownarrows \\ \Vc &\rightrightarrows &\Pc
\end{matrix}$ be a double groupoid. Let $\Box(\Vc, \Hc)$ be the
coarse double groupoid with sides in $\Vc$ and $\Hc$ and let $\Pi:
\B \to \Box(\Vc, \Hc)$ be the natural morphism of double
groupoids. Let also
\begin{align*} \Kc & : = \{K \in
\B: \; t(K), b(K),l(K), r(K)  \in \Pc \},
\\ \F &: = \text{image of } \Pi;
\end{align*}
$\F$ is called the frame of $\B$. The frame is what we call a
\emph{slim} double groupoid; that is, a double groupoid in which
every box is uniquely determined by its edges. On the other hand,
$\Kc$ is an abelian group bundle over $\Pc$ with respect to
vertical composition, which coincides with the horizontal one in
$\Kc$.

\begin{mtheorem}\label{teo:main}(a) \,
 The structure of the double groupoid $\B$ is determined
completely by its associated abelian group bundle $\Kc$ and  its
frame $\F$, together with some extra cohomological data.

\medbreak (b) \, Assume that $\B$ is slim and satisfies the
filling condition. Then there exists a groupoid $\D$ and morphisms
of groupoids $i: \Hc \to \D$, $j: \Vc \to \D$ with $\D =
j(\Vc)i(\Hc)$ (unique up to isomorphisms), such that $\B \simeq
\square(\D, j, i)$-- as defined in Subsection \ref{doubleoffact}.
\qed \end{mtheorem}

Part (a) of Theorem \ref{teo:main} is contained in Theorem
\ref{main1}. Part (b) is contained in Theorem \ref{class-slim}.

\medbreak R. Brown has kindly pointed out to us that the
construction of the double groupoid $\square(\D, j, i)$ was found
by him a long time ago, and that it was taken up by Lu and
Weinstein \cite{wl}.

\medbreak It is well-known that a discrete groupoid can be fully
described in group-theoretical terms. Our main result shows that
there is an analogous description, albeit more complicated, of the
structure of discrete double groupoids in group-theoretical terms.
Part (a) of Theorem \ref{teo:main}, included to underline the
completeness of our approach, does not require the filling
condition and suggests a study of the corresponding cohomological
theory, that we postpone to a future publication. Part (b) of
Theorem \ref{teo:main} gives a description of a much larger class
of double groupoids than the vacant double groupoids characterized
in \cite[2.15]{mk1}. Namely, a vacant double groupoid is always
slim, and if $\B$ is a slim double groupoid, then the following
are equivalent:
\begin{enumerate}
\item[(i)] $\B$ is vacant.

\item[(ii)] $\B$ satisfies the filling condition, the morphisms $i$
and $j$ are injective, and the factorization is exact.
\end{enumerate}

Thus, there are plenty of slim double groupoids satisfying the
filling condition that are not vacant.

\medbreak The extension of our main result to the context of Lie
double groupoids, or other natural settings, requires some extra
labor. The definition of the cocycles in part (a) requires the
choice of a section, imposing an extra obstruction. The study of
part (b) in the context of Lie double groupoids, under suitable
natural hypotheses, has been recently carried out in \cite{aot},
starting from our main result in Theorem \ref{class-slim}.

\medbreak In the paper \cite{AN1} we proved that  vacant finite
double groupoids gave rise, in a natural way, to a class of tensor
categories.  Thus the results of \cite{AN1} generalized a
well-known construction in Hopf algebra theory studied, among
others, by G. I. Kac, Majid and Takeuchi.

Later, in \cite{AN2}, this result  was extended to the much more
general class of finite double groupoids satisfying only the
filling condition \ref{filling}. It turns out that double
groupoids giving rise through this construction to a special class
of tensor categories called \emph{fusion} categories must be slim.
We discuss this class of double groupoids in more detail in the
last section of the paper. We plan to apply the results of this
last section to the determination of the corresponding fusion
categories in a subsequent publication.

\medbreak This paper is organized as follows. In Section
\ref{do-gpd}, we recall the definition and special features of
double groupoids.  In Section \ref{fusion} we discuss a special
class of slim double groupoids motivated by the paper \cite{AN2}.
In the Appendix we discuss an alternative approach to the corner
functions introduced in \cite{AN2}.

\section*{Acknowledgements}

Several of the results in this paper were obtained while the
authors visited the Institut des Hautes \' Etudes Scientifiques
(IHES) in France. Both authors are grateful to the Institute for
the excellent atmosphere. Part of the work of the first author was
done while visiting Hunan University (China); he thanks Prof.
Zhang Shouchuan for his warm hospitality. Part of the work of the
second author was done while visiting Nanjing Agricultural
University (China); she thanks Prof. Zhang Lianyun for his warm
hospitality. Both authors thanks Prof. Ronnie Brown for
interesting discussions and many bibliographical references.

\subsection*{Notation}
The cardinality of a set $S$ will be indicated by $|S|$.

Throughout this paper we fix a nonempty set $\Pc$. A groupoid $\G$
on the base $\Pc$, with source $s$ and target $e$, will be
indicated by $s, e: \G \rightrightarrows \Pc$ or simply $\G$ if no
confusion arises. Composition in $\G$ will be indicated by
juxtaposition of arrows; so that if $g, h \in \G$, such that $e(g)
= s(h)$, their composition will be denoted by $gh \in \G$. Let $P,
Q\in \Pc$. As usual, $\G(P, Q)$ is the set of arrows from $P$ to
$Q$ and $\G(Q) = \G(Q, Q)$. A map $p:\Ee \to \Pc$ will be called a
fiber bundle. Let $p: \Ee \to \Pc$ be a fiber bundle. Recall that
a \emph{left action} of $\G$ on $p$ is a map $\fde : \G
\Times{e}{p} \Ee \to \Ee$ such that
$$ p(g\fde x) = s(g), \qquad g \fde(h \fde x)  = gh \fde x,\qquad
 \iddh \,{p(x)} \, \fde  x  =  x,
$$
for all $g, h \in \G$, $x \in \Ee$  composable in the appropriate
sense.  Hence, if $\Ee_b := p^{-1} (b)$, then the action of $g\in
\G$ is an isomorphism $g\fde \underline{\quad}: \Ee_{t(g)} \to
\Ee_{s(g)}$. Let $x\in \Ee$ and define
\begin{align*}
\Oc_x &= \{g \fde x: g\in \G, e(g) = p(x)\}, & \qquad &\text{the
orbit of } x, \\
\G^x &= \{g\in \G: g \fde x = x\} < \G(p(x)), & \qquad &\text{the
isotropy subgroup of } x.\end{align*}

\section{Double groupoids}\label{do-gpd}

Let $\T$ be a  \emph{double groupoid} \cite{ehr}; we follow the
conventions and notations from \cite[Section 2]{AN1} and
\cite[Section 1]{AN2}.  As usual, we represent $\T$ in the form of
four related groupoids
$$\begin{matrix} \B &\rightrightarrows &\Hc
\\\downdownarrows &&\downdownarrows \\ \Vc &\rightrightarrows &\Pc \end{matrix}$$
subject to a set of axioms. The source and target maps  of these
groupoids are indicated by $\prin, \fin: \B \to \Hc$; \, $\pri,
\fine: \B \to \Vc$; \, $\pri, \fine: \Hc \to \Pc$; \, $\prin,
\fin: \Vc \to \Pc$ (`top', `bottom', `right' and `left'). An
element $A\in \B$ is depicted as a box $$A =
\begin{matrix} \quad t \quad \\ l \,\, \boxe \,\, r \\ \quad b\quad
\end{matrix}$$ where $t(A) = t$, $b(A) = b$, $r(A) = r$, $l(A) = l$,
and the four vertices of the square representing $A$ are $tl(A) =
lt(A)$, $tr(A) = rt(A)$, $bl(A) = lb(A)$, $br(A) = rb(A)$. The notation
$A \vert B$ means that $r(A)=l(B)$ ($A$ and $B$ are
horizontally composable); the corresponding horizontal product is denoted
$AB$. Similarly, $\displaystyle\frac{A}{B}$ means that
$b(A)=t(B)$ ($A$ and $B$ are vertically composable)
and the vertical product is denoted
$\begin{matrix} A\\B \end{matrix}$.

The notation
$A = $ \begin{tabular}{|p{0,1cm}|} \hhline{|=|} \\
\hline\end{tabular}\, means that $t(A)$ is an identity;
analogously, $B=$ \begin{tabular}{||p{0,1cm}|} \hline \\
\hline \end{tabular} means that $l(B)$ is an identity, etc.

These four groupoids should satisfy certain axioms, see e. g.
\cite{AN1}. In particular, $\idd \, {\id_\Hc P} = \idd \, {\id_\Vc
P}$, for any $P\in \Pc$; this box is denoted $\Theta_P$ and
clearly it is of the form \begin{tabular}{||p{0,1cm}||} \hhline{|=|} \\
\hhline{|=|}\end{tabular}.

\subsection{Core groupoids}\label{core}

\

We recall the core groupoid  $\Ec$ of $\T$ introduced by
Brown and Mackenzie \cite{mk2, BM}. See also \cite{AN2}. Let
\begin{align*} \Ec & : = \{ E \in
\B: \; r(E), t(E) \in \Pc \}.\end{align*}

Thus, elements of  $\Ec$ are of the form
\begin{tabular}{|p{0,1cm}||} \hhline{|=||} \\
\hline \end{tabular}. There is a  groupoid structure $s, e: \Ec
\rightrightarrows \Pc$,  $s(E) = bl(E)$, $e(E) = br(E)$, $E \in
\Ec$, identity map $\id: \Pc \to \Ec$, $P \mapsto \Theta_P$, and
composition $\Ec \Times{e}{s} \Ec \to \Ec$,  given by
\begin{equation}\label{productcore} E \circ M : =
\left\{\begin{matrix}\iddv l(M)& M \vspace{-4pt}\\ E &\iddv
b(M)\end{matrix} \right\},
\end{equation} $M, E \in \Ec$.
The inverse of  $E \in \Ec$ is $E^{(-1)}: = (E \iddv b(E)^{-1})^v
= \left\{\begin{matrix}\iddv l(E)^{-1} \vspace{-4pt}\\
E^h\end{matrix} \right\}$. If $Q\in \Pc$, the group $\Ec(Q)$
consists of boxes in $\B$ of the form
$\begin{matrix} h \hspace{-0.5cm} & \begin{tabular}{|p{0,1cm}||} \hhline{|=||} \\
\hline \end{tabular} \\  & y\end{matrix}$, with $y \in \Hc(Q)$,
$h\in \Vc(Q)$.

\medbreak If $(x,g)\in \Hc \Times{r}{t}\Vc$ and $B\in \B$ then we
consider the sets of boxes with fixed `upper and right' sides
\begin{equation}\label{ur}
\Ur(x,g) = \Big\{ U \in \B: \, U =
\begin{matrix} \quad x \,\,\, \quad \\  \,\, \boxe \,\, g \\ \quad \quad \end{matrix}
\Big\}, \quad \Ur(B) = \Ur(t(B), r(B)). \end{equation}

Let $\gamma: \B \to \Pc$ be the `left-bottom' vertex, $\gamma(B) =
lb(B)$.

\begin{prop}\label{prop-core} (a). There is an action of  $\Ec$ on $\gamma: \B \to \Pc$
given by \begin{equation}\label{actioncore} E \acts A : =
\left\{\begin{matrix}\iddv l(A)& A \vspace{-4pt}\\ E &\iddv
b(A)\end{matrix} \right\}, \quad A\in \B, E \in \Ec.
\end{equation}

\medbreak (b). Let $B\in \B$. Then  $\Oc_B = \Ur(B)$ and $\Ec^B$
is trivial.
\end{prop}

\pf (a) is straightforward. (b): it follows from Definition
\eqref{actioncore} that $\Oc_B \subseteq \Ur(B)$. Then observe
that for any $C\in \Ur(B)$, there exists a unique $E\in \Ec$ such
that $E \acts B = C$, namely $E = \left\{\begin{matrix}B^v
\vspace{-4pt}\\ C\end{matrix} \right\} \id b(B)^{-1}$. \epf

\medbreak
\subsection{The associated abelian group bundle}\label{pith}

\

An important invariant of a double groupoid is the intersection
$\Kc$ of all four core groupoids:
\begin{align*} \Kc & : = \{K \in
\B: \; t(K), b(K),l(K), r(K)  \in \Pc \}.\end{align*}
Thus a box is in $\Kc$ if and only if it is  of the form \begin{tabular}{||p{0,1cm}||} \hhline{|=|} \\
\hhline{|=|}\end{tabular}. Let $p: \Kc \to \Pc$ be the `common
vertex' function, say $p(K) = lb(K)$. For any $P\in \Kc$, let
$\Kc(P)$ be the fiber  at $P$; $\Kc(P)$ is an abelian group under
vertical composition, that coincides with horizontal composition. This
is just the well-known fact: ``a double group is the same as an
abelian group". Indeed, apply the interchange law
\begin{equation*}
\begin{matrix} (K L) \vspace{-2pt}\\ (M N) \end{matrix} =
\left(\begin{matrix} K  \vspace{-4pt}\\ M   \end{matrix}\right)
\left(\begin{matrix} L \vspace{-4pt}\\ N  \end{matrix}\right)
\end{equation*}
to four boxes $K,L,M,N \in \Kc(P)$: if $L = M = \Theta_P$, this
says that $\begin{matrix} K  \vspace{-4pt}\\ N \end{matrix}  = KN$
and the two operations coincide.  If, instead, $K = N = \Theta_P$,
 this says that $\begin{matrix} L  \vspace{-4pt}\\ M \end{matrix}
= ML$, hence the composition is abelian. Note that this operation
in $\Kc$ coincides also with the core multiplication
\eqref{productcore}. In short, $\Kc$ is an abelian group bundle
over $\Pc$.

\medbreak The vertical and the horizontal groupoids $\Vc$ and
$\Hc$ act on $\Kc$ by vertical, respectively horizontal,
conjugation:
\begin{flalign}
\label{accionver} & \text{If } g\in \Vc(Q,P) \text{ and } K\in
\Kc(P) \text{ then } g\cdot K :=
  \begin{matrix} \id g\vspace{-1pt}\\K  \vspace{-4pt}\\ \id g^{-1} \end{matrix}
\in \Kc(Q).&
\\
\label{accionhor} &\text{If } x\in \Hc(Q,P) \text{ and } K\in
\Kc(P) \text{ then } x\cdot K =  \id xK\id x^{-1}\in \Kc(Q).&
\end{flalign}
Both actions are by group bundle automorphisms.

\subsection{Frame of a double groupoid}

\

The frame of a double groupoid, as said in the Introduction, is
the image of the natural map into the coarse double groupoid of
the side groupoids. In this subsection, we made this statement
precise; then we show that the initial double groupoid is
determined by its frame, its associated abelian group bundle and
some actions and cocycles. In the next subsection we explain this
fact in terms of extensions.

Let $\Pc$ be a set and $\Vc$, $\Hc$ be groupoids over $\Pc$
denoted vertically and horizontally, respectively. Let $\Box(\Vc,
\Hc)$ be the set of quadruples $\begin{pmatrix} \quad x  \quad \\  f \quad g \\
\quad y \quad\end{pmatrix}$ with $x,y\in \Hc$, $f,g\in \Vc$ such
that
\begin{equation*}
l(x) = t(f), \quad r(x) = t(g), \quad l(y) = b(f), \quad r(y) =
b(g). \end{equation*} If no confusion arises, we shall denote a
quadruple as above by  a box $\begin{matrix} \quad x \quad \\
h\,\, \boxe \,\, g \\ \quad y \quad
\end{matrix}$.

The collection $\begin{matrix} \Box(\Vc, \Hc) &\rightrightarrows
&\Hc \\\downdownarrows &&\downdownarrows \\ \Vc &\rightrightarrows
&\Pc \end{matrix}$ forms a double groupoid, called the
\emph{coarse double groupoid} with sides in $\Vc$ and $\Hc$, with
horizontal and vertical compositions given by
\begin{equation*} \begin{matrix} \quad x \quad
\\  h\,\, \boxe  \,\, g \\ \quad y \quad
\end{matrix} \,  \begin{matrix} \quad x' \quad
\\  g\,\, \boxe  \,\, g' \\ \quad y' \quad
\end{matrix}  = \begin{matrix} \quad xx' \quad
\\  h\,\, \boxe  \,\, g' \\ \quad yy' \quad
\end{matrix}, \quad \begin{matrix} \quad x \quad
\\  h\,\, \boxe  \,\, g \\ \quad y \quad \\  \quad y\quad \\ h'\,\, \boxe  \,\, g' \\ \quad y' \quad \end{matrix}  =
\begin{matrix} \quad x \quad
\\  hh'\,\, \boxe  \,\, gg' \\ \quad y', \qquad
\end{matrix},\end{equation*} for all $x, y, x', y' \in \Hc$, $g, h, g', h' \in \Vc$ appropriately composable.

Let $\begin{matrix} \B &\rightrightarrows &\Hc
\\\downdownarrows &&\downdownarrows \\ \Vc &\rightrightarrows &\Pc \end{matrix}$ be a double
groupoid. There is a map $\proy: \B \to \Box(\Vc, \Hc)$ given by
$$
\proy\left(\begin{matrix} \quad x \quad \\ f \,\, \boxe \,\, g \\
\quad y \quad \end{matrix} \right) = \begin{pmatrix} \quad x \quad
\\  f \quad g \\ \quad y \quad\end{pmatrix}, \qquad
 \begin{matrix}
\quad x \quad \\ f \,\, \boxe \,\, g \\ \quad y \quad \end{matrix}
\in \B.
$$

\medbreak Clearly, $\proy$ induces a morphism of double groupoids
$\B \to  \Box(\Vc, \Hc)$.

\begin{definition}We shall say that $\B$ is \emph{slim} if $\proy$ is injective
(any box is determined by its four sides). \end{definition}

\medbreak Let $\F$ be the image of $\proy$. The \emph{frame} of
$\T$ is the slim double groupoid $$\begin{matrix} \F
&\rightrightarrows &\Hc
\\\downdownarrows &&\downdownarrows \\ \Vc &\rightrightarrows &\Pc
\end{matrix}.$$

Let $\B$ be a double groupoid. Several properties of $\T$ are
controlled by its frame $\F$. Recall the following definitions
\cite[Definition 2.3]{BM}.
\begin{enumerate}
\item[(a)] $\T$ is \emph{horizontally transitive} if every configuration of matching sides
$\boxu$ can be completed to at least one box in $\B$.

\item[(b)] $\T$ is \emph{vertically transitive} if every configuration of matching sides
$ \boxa $ can be completed to at least one box in $\B$.

\item[(c)] $\T$ is \emph{transitive} or \emph{locally trivial}
if it is both vertically  and horizontally transitive.
\end{enumerate}

\begin{obs} Let $\T$ be a double groupoid. Then

\indent (i) $\T$ satisfies the filling condition \ref{filling} if
and only if $\F$ does so.

\indent (ii) $\T$ is horizontally (vertically) transitive if and
only if so $\F$ is so.

\indent (iii) If $\T$ is vacant then it is slim.
 \end{obs}

We are now ready to prove our first basic result. \emph{Let us fix
a section $\mu: \F \to \B$ of $\proy$}. Recall the action $\acts$
of the core groupoid $\Ec$ on $\gamma: \B \to \Pc$ defined in
\eqref{actioncore}.

\begin{lema}\label{framepith} The map
$\Psi:\Kc\Times{p}{\gamma} \F \to \B$ given by $\Psi(K, F) =
K\acts \mu(F)$ is a bijection.
\end{lema}

\pf In other words, we have to show that for any $B \in
\Ur(\mu(F))$, there is a unique $K \in \Kc$ such that $B = K\acts
\mu(F)$. Note that $\proy (K\acts \mu(F)) = F$. Hence, $K\acts
\mu(F) = K'\acts \mu(F')$ implies $F = F'$; thus $K\acts \mu(F) =
K'\acts \mu(F)$, and $K = K'$ by Proposition \ref{prop-core} (b).
That is, $\Psi$ is injective. We show that it is surjective. Let
$B\in \B$ and let $F = \proy(B)$. Since $B$ and $\mu(F)$ have the
same sides, there exists $K\in \Ec$ such that $B = K\acts \mu(F)$,
again by Proposition \ref{prop-core}. But clearly $K\in \Kc$. \epf

We next introduce vertical and horizontal cocycles to control the
lack of multiplicativity of the section $\mu$.  We define $\tau:
\F \Times{r}{l}\F \to \Kc$ and $\sigma: \F \Times{b}{t}\F \to \Kc$
by
\begin{align}\label{deftau}
\mu(F)\mu(G)&=\tau(F, G) \acts \mu(F G), \qquad r(F) = l(G),
\\\label{defsigma}
\begin{matrix} \mu(F)  \vspace{-4pt}\\ \mu(G) \end{matrix} & =
\sigma(F, G) \acts \mu\left(\begin{matrix} F\vspace{-4pt}\\ G
\end{matrix}\right), \qquad b(F) = t(G).
\end{align}
That is,
\begin{align*}
\mu(F)\mu(G)=\left\{\begin{matrix}\id l(F) & \mu(F G) \\
\tau(F, G) & \id b(FG) \end{matrix}\right\}, \qquad
\begin{matrix} \mu(F)  \vspace{-4pt}\\ \mu(G) \end{matrix}  = \,\left\{
\begin{matrix}\id l(F)l(G) & \mu\left(\begin{matrix} F\vspace{-4pt}\\ G \end{matrix}\right)
\\\sigma(F, G) & \id b(G)\end{matrix}\right\},
\end{align*}
for appropriate $F$, $G\in \F$. The cocycles $\sigma$ and $\tau$
are well-defined in virtue of Lemma \ref{framepith}. If we assume
that $\mu(\id x) = \id x$ and $\mu(\id g) = \id g$ for any $x\in
\Hc$ and $g\in \Vc$ then $\sigma$ and $\tau$ are normalized:
\begin{align}\label{tauunitario}
\tau(F, \id r(F)) &= \Theta_{bl(F)}= \tau(\id l(F), F),
\\\label{sigmaunitario}
\sigma(F, \id b(F)) &= \Theta_{bl(F)}= \sigma(\id t(F), F).
\end{align}

Now we can reconstruct the horizontal and vertical products of
$\B$ in terms of the associated abelian group bundle $\Kc$, the
frame slim double groupoid $\F$, the actions \eqref{accionhor},
\eqref{accionver} and the cocycles $\sigma$ and $\tau$. If $K,L\in
\Kc$, $F,G\in \F$ then

\begin{align}\label{prodhorreconstr}
\big(K\rightharpoondown \mu(F) \big) \,
\big(L\rightharpoondown \mu(G) \big)  &=
\big(K (b(F)\cdot L) \tau(F, G) \big)
\rightharpoondown \mu(FG),
\end{align}
if $r(F) = l(G)$ and
\begin{align}
\label{prodverreconstr} \
\left\{\begin{matrix}K\rightharpoondown \mu(F) \\
L\rightharpoondown \mu(G)\end{matrix}\right\} &=
\big((l(g)^{-1}\cdot K)\, L\, \sigma(F, G) \big)
\rightharpoondown \mu\left(\begin{matrix}F \\ G\end{matrix}\right)
\end{align}
if $b(F) = t(G)$.

\subsection{Extensions of double groupoids by abelian group bundles}

\

Our aim now is to interpret Lemma \ref{framepith} and formulas
\eqref{prodhorreconstr}, \eqref{prodverreconstr} in cohomological
terms. The description in the preceding subsection suggests the
following construction.

Let $\begin{matrix} \F &\rightrightarrows &\Hc
\\\downdownarrows &&\downdownarrows \\ \Vc &\rightrightarrows &\Pc
\end{matrix}$ be any double groupoid (not necessarily
slim) and let $\gamma:\Kc \to \Pc$ be any abelian group bundle.
Assume that $\Vc$ and $\Hc$ act on $\Kc$ by group bundle
isomorphisms. Let $\tau: \F \Times{r}{l}\F \to \Kc$ and $\sigma:
\F \Times{b}{t}\F \to \Kc$ be maps such that
\begin{align}\label{sigma-inicial}
\gamma(\sigma(F,G)) &= bl(G), \qquad \text{if } b(F) = t(G),
\\\label{tau-inicial}
\gamma(\tau(F, G)) &= bl(F), \qquad \text{if } r(F) = l(G),
\end{align}
normalized by \eqref{tauunitario} and \eqref{sigmaunitario}.
Consider the collection $\begin{matrix}
\Kc\Times{p}{\gamma} \F &\rightrightarrows &\Hc
\\\downdownarrows &&\downdownarrows \\ \Vc &\rightrightarrows &\Pc
\end{matrix}$, where:

\begin{itemize}
  \item The maps $t,b,l,r$ on $\Kc\Times{p}{\gamma} \F$ are
  defined by those in $\F$: $t(K, F) = t(F)$ and so on.

\item The horizontal and vertical products in
  $\Kc\Times{p}{\gamma} \F$ are given by
\begin{align}\label{prodhorreconstrbis}
(K,F)(L,G) &= \big(K \,(b(F)\cdot L)\, \tau(F, G), FG\big),
\qquad \text{if } F\vert G, \\
\label{prodverreconstrbis}
\begin{matrix}(K,F) \\(L,G) \end{matrix} &=
\left( (l(G)^{-1}\cdot K)\, L\, \sigma(F, G),
\begin{matrix}F \\ G\end{matrix}\right), \qquad \text{if } \dfrac FG.
\end{align}

\item The identity maps $\id: \Vc\to\Kc\Times{p}{\gamma} \F$, $\id: \Hc\to\Kc\Times{p}{\gamma} \F$ are
  given by $\id g = (\Theta_{b (g)}, \id g)$, $\id x = (\Theta_{l (x)}, \id
  x)$, $g\in \Vc$, $x\in \Hc$.

\item The inverse of $(K, F)$ with respect to the horizontal and
vertical products are respectively given by
\begin{align}\label{inversa-hor}
(K,F)^h &= \left(b(F)^{-1}\cdot\left( K^{-1} \tau(F, F^h)^{-1}\right), F^h\right),  \\
\label{inversa-ver} (K,F)^v &= \left(\left(l(F)\cdot
K\right)^{-1}\sigma(F, F^v)^{-1}, F^v\right).
\end{align}
\end{itemize}

\begin{prop}\label{prop-extension} $\Kc\Times{p}{\gamma} \F$ is a
double groupoid if and only if, for all $F,G,H\in \F$,
\begin{align}\label{taucocyclo}
\tau(F, G)\,\tau(FG, H)&=\tau(F, GH) \, \big( b(F) \cdot \tau(G, H)\big),
\quad F\vert G \vert H;
\\\label{sigmacocyclo}
\sigma(G, H) \, \sigma\left(F, \begin{matrix} G \vspace{-4pt}\\ H
\end{matrix}\right) & = \big( l(H)^{-1} \cdot\sigma(F, G)\big) \,  \sigma\left(\begin{matrix}
F \vspace{-4pt}\\ G \end{matrix}, H\right), \quad \begin{tabular}{p{0,4cm}}$F$ \\ \hline $G$\\
\hline $H$\end{tabular};
\\\label{compatib-acciones} l(H)^{-1} \cdot (t(H)\cdot L) &= b(H) \cdot (r(H)^{-1}\cdot
L), \quad L\in \Kc(tr(H));
\end{align}
\begin{multline}
\label{sigma-tau-compatibles} \big(l(H)^{-1}\cdot \tau(F,G)\big)
\tau(H,J)\sigma(FG, HJ) \\= \big(b(H)\cdot \sigma(G,J)\big)
 \sigma(F,H)
\tau\left(\begin{matrix} F \vspace{-4pt}\\ H
\end{matrix}, \begin{matrix} G\vspace{-4pt} \\ J \end{matrix} \right), \quad
\begin{tabular}{p{0,4cm}|p{0,4cm}} $F$ & $G$ \\ \hline $H$ & $J$
\end{tabular}.
\end{multline}
\end{prop}

\begin{definition} If the conditions in Proposition \ref{prop-extension} hold, we say that the double groupoid
$\Kc\Times{p}{\gamma} \F$ is an \emph{abelian extension} of the
abelian group bundle $\Kc$ by $\F$.
\end{definition}

\pf The associativity of the horizontal and vertical compositions
are respectively equivalent to \eqref{taucocyclo} and
\eqref{sigmacocyclo}.

We have to check the axioms of double groupoid as in \cite{bs}; we
follow  \cite[Lemma 1.2]{AN1}. All the axioms are consequences of
the definitions (since the axioms hold in $\F$) except the
interchange law, which is equivalent to \eqref{compatib-acciones}
and \eqref{sigma-tau-compatibles}. Indeed, let $H\in \F$ and $L\in
\Kc(tr(H))$. Computing
$\left\{\begin{matrix} \id t(H) & L \\ H & \id
r(H) \end{matrix}\right\}$ in two different ways, we see that
\eqref{compatib-acciones} is equivalent to the interchange law in
this case.
Next, consider

\centerline{$(K,F)$, $(L,G)$, $(M,H)$, $(N,J)\in
\Kc\Times{p}{\gamma} \F$ such that
$\begin{tabular}{p{0,4cm}|p{0,4cm}} $F$ & $G$ \\
\hline $H$ & $J$ \end{tabular}$.}

Compute $\begin{matrix}
\{(K,F)(L,G)\} \vspace{-4pt} \\ \{(M,H)(N,J)\} \end{matrix}$ and
$\left\{\begin{matrix} (K,F)  \vspace{-6pt} \\
(M,H)\end{matrix}\right\} \left\{\begin{matrix} (L,G) \vspace{-6pt} \\
(N,J) \end{matrix}\right\}$. The resulting
expressions are equal if and only
if
\begin{multline*}
l(H)^{-1} \cdot (b(F)\cdot L) l(H)^{-1}\cdot \tau(F,G)
\tau(H,J)\sigma(FG, HJ) \\ = b(H) \cdot (l(J)^{-1}\cdot L)
b(H)\cdot \sigma(G,J) \sigma(F,H) \tau\left(\begin{matrix} F \\ H
\end{matrix}, \begin{matrix} G \\ J \end{matrix} \right)
\end{multline*}
It is not difficult to see that this is equivalent to
\eqref{compatib-acciones} and \eqref{sigma-tau-compatibles}. \epf

\begin{obs} As we shall see later, see Section \ref{slim-gpd},
the slim double groupoid $\F$ determines a 'diagonal' groupoid
$\D$ endowed with groupoid maps $j: \Vc \to \D$, $i: \Hc \to \D$.
In terms of this groupoid, Condition \eqref{compatib-acciones}
means that the actions of $\Hc$ and $\Vc$ come from an action of
$\D$ on $\Kc$.
\end{obs}

\begin{obs} Conditions \eqref{taucocyclo} and \eqref{sigmacocyclo}
are cocycle conditions on the horizontal and vertical composition
groupoids. Together with \eqref{sigma-tau-compatibles} they give a
cocycle condition in the double complex associated to the double
groupoid $\F$ as considered in \cite{AN1, AM}. \end{obs}

Assume that the hypotheses in Proposition \ref{prop-extension} are
fulfilled. Let us identify $\Kc$ with a subset of
$\Kc\Times{p}{\gamma} \F$ via $K\mapsto (K, \Theta_P)$, if $K\in
\Kc(P)$. Also, let $\mu: \F \to \Kc\Times{p}{\gamma} \F$, $\mu(F)
= (\Theta_{bl(F)}, F)$. Then
$$
(K, F) = \begin{matrix}\id l(F) & \mu(F) \\
K & \id b(F) \end{matrix}, \qquad \text{for any } (K, F)\in
\Kc\Times{p}{\gamma} \F.
$$
Hence the formulas \eqref{prodhorreconstrbis} and
\eqref{prodverreconstrbis} are equivalent to
\eqref{prodhorreconstr} and \eqref{prodverreconstr}, respectively.
In particular we have

\begin{theorem}\label{main1} Any double groupoid is an abelian extension of
its associated abelian group bundle by its frame. \qed
\end{theorem}

\section{Slim double groupoids and factorizations of
groupoids}\label{slim-gpd}

Let $\Vc$ and $\Hc$ be groupoids over $\Pc$. In this section we
describe all slim double groupoids satisfying the filling
condition \ref{filling} whose groupoids of vertical and horizontal
edges coincide with $\Vc$ and $\Hc$, respectively.

\subsection{Double groupoid associated to a diagram of groupoids}\label{doubleoffact}

\

Let us say that a \emph{diagram} $(\D, j, i)$ over $\Vc$ and $\Hc$
is a groupoid $\D$ over $\Pc$ endowed with groupoid maps over
$\Pc$
$$i: \Hc \to \D, \quad j: \Vc \to \D.$$
 The class of all diagrams $(\D, j, i)$ over $\Vc$ and $\Hc$
is a category with morphisms $(\D, j, i) \to (\D', j', i')$ being
morphisms $f: \D \to \D'$ of groupoids over $\Pc$ such that $fi =
i'$ and $fj = j'$.

Consider the full subcategory of diagrams $(\D, j, i)$ with $\D =
j(\Vc)i(\Hc)$; that is, such that every arrow in $\D$ can be
written as a product $j(g)i(x)$, for some $g \in \Vc$, $x \in
\Hc$, with $b(g) = l(x)$. An object in this subcategory will be
called a $(\Vc, \Hc)$-\emph{factorization} of $\D$.

\medbreak Each diagram $(\D, j, i)$ has an associated double
groupoid $\square(\D, j, i)$ defined as follows. Boxes in
$\square(\D, j, i)$ are of the form
$$A =
\begin{matrix} \quad x \quad \\ h \,\, \boxe \,\, g \\ \quad y\quad
\end{matrix} \in \square(\Vc, \Hc),$$ with $x, y \in \Hc$, $g, h \in \Vc$, such that $$i(x)j(g) = j(h) i(y)\quad \text{ in } \D.$$

Notice that $\square(\D, j, i)$ is stable under vertical and
horizontal products in $\square(\Vc, \Hc)$; therefore it is itself
a double groupoid. By its very definition, $\square(\D, j, i)$ is
\emph{slim}.

\medbreak The assignment $(\D, j, i) \to \square(\D, j, i)$ just
defined is clearly functorial.

\begin{exa}Let $G$ be a simply connected Poisson-Lie group,
$\mathfrak g$ its Lie algebra. It is well known that $\mathfrak g$
is a  Lie bialgebra; let $\mathfrak g^*$ be the dual Lie algebra
and let $\mathfrak d$ be the correspondig Drinfeld double. Let
$G^*$ and $D$ be simply connected Lie groups with Lie algebras
$\mathfrak g^*$ and $\mathfrak d$, respectively. Then the maps $G
\to D$ and $G^* \to D$ give rise to a double symplectic groupoid
\cite[Theorem 3]{wl}.
\end{exa}

\begin{lema} The core groupoid of $\square(\D, j, i)$ is isomorphic to the groupoid $\Vc^{\op} {}_{j}\times_{i}\Hc: = \{(g, x) \in \Vc^{\op} {}_{b}\times_{l}\Hc: \, j(g) = i(x^{-1})\} \subseteq \Vc^{\op} {}_{b}\times_{l}\Hc$. \end{lema}

\pf  An isomorphism is given by the map $\Ec \to \Vc^{\op}
{}_{j}\times_{i}\Hc$, defined by $E \mapsto (l(E), b(E))$. This
map is surjective by construction of $\square(\D, j, i)$; it is
injective as a consequence of the slim condition on $\square(\D,
j, i)$.  \epf

\begin{obs} If $\D = j(\Vc)i(\Hc)$ is  a factorization, then
$\square(\D, j, i)$ satisfies the filling condition \ref{filling}.
Indeed, if $g \in \Vc$, $x \in \Hc$, are such that $r(x) = t(g)$,
then the condition $\D = j(\Vc)i(\Hc)$ implies that there exist $y
\in \Hc$, $h \in \Vc$, such that $j(h)i(y) = i(x)j(g)$. Then, by
construction, the box $\begin{matrix} \quad x \quad
\\ h \,\, \boxe \,\, g \\ \quad y\quad
\end{matrix}$ is a filling in $\square(\D, j, i)$ for $\begin{matrix}\;  x
\quad \\ \boxur \; g \\ \quad \end{matrix}$.
\end{obs}

\begin{exa} Suppose $\D = \Pc^2$ is the coarse groupoid on $\Pc$.
Let the maps $i: \Hc \to \D$, $j: \Vc \to \D$, be defined by $i(x)
= (l(x), r(x))$, $x\in \Hc$, and $j(g) = (t(g), b(g))$, $g\in \Vc$.
Let $\B = \square(\Pc^2, j, i)$ be the
associated double groupoid. The relations $i(x)j(g) = j(h)i(y)$
are satisfied in $\D$, for all $x, y \in \Hc$, $g, h \in \Vc$,
such that $r(x) = t(g)$, $b(h) = l(y)$, $l(x) = t(h)$, $r(y) =
b(g)$.  Hence, for all such $x, y, g, h$ there is a box
$\begin{matrix} \quad x \quad
\\ h \,\, \boxe \,\, g \\ \quad y\quad
\end{matrix}$ in $\B$. According to the composition rules in $\B$,
it turns out that $\B$ is exactly the \emph{coarse} double
groupoid $\square(\Vc, \Hc)$. \end{exa}

\medbreak We shall show that  double groupoids of the form
$\square(\D, j, i)$ exhaust the class of slim double groupoids
which satisfy the filling condition \ref{filling}.

\subsection{Free product of groupoids}\label{presentation}

\

We now need to recall the basic properties of the free product
construction for groupoids. We refer the reader to \cite{higgins}
for a detailed exposition; this will be used in the definition of
the diagonal groupoid in Subsection \ref{diagonalofslim} and in
the proof of Theorem \ref{mainresult}.

In this subsection we work in the category of groupoids over
$\Pc$; morphisms in this category are the identity on $\Pc$.
Subgroupoids with the same base $\Pc$ are often called `wide'.

Let $\Vc$ and $\Hc$ be groupoids, alluded to as `vertical' and
`horizontal', respectively. Let $\Vc = \langle X | R \rangle$,
$\Hc = \langle Y| S \rangle$, be presentations of $\Vc$ and $\Hc$
by generators and relations \cite[Chapter 9]{higgins}. Let $\Vc
* \Hc = \langle X \coprod Y | R \coprod S \rangle$ be the
\emph{free product }of the groupoids $\Vc$ and $\Hc$; $\Vc * \Hc$
is the coproduct of $\Vc$ and $\Hc$ in the category of groupoids
over $\Pc$. In other words, the groupoid $\Vc * \Hc$ is
characterized by the following universal property: for every
groupoid $\G$ and groupoid maps $i: \Hc \to \G$, $j: \Vc \to \G$,
there is a unique morphism of groupoids $f: \Vc * \Hc \to \G$ such
that $f\vert_{\Vc} = j$, and $f\vert_{\Hc} = j$. In particular it
does not depend on the choice of the presentations of $\Vc$ and
$\Hc$.

Note that our free product of $\Vc$ and $\Hc$ is close to, but not
the same as, the free product in \cite[Chapter 9]{higgins};
precisely, it is the free product with amalgamation of identities
from {\it loc. cit.}

\medbreak An alternative way of describing the free product is the
following. Consider the set $\path(\Q)$ of all paths of the quiver
$\Q = \Hc \coprod \Vc$. An element in $\path(\Q)$ is either a an
element $P \in \Pc$ that will be indicated by $[P]$, or a sequence
$(u_1, \dots, u_n)$, $n \geq 1$, with $u_i \in \Q$, $e(u_i) =
s(u_{i+1})$.

A path $U \in \path(\Q)$ is called \emph{reduced} if either $U =
[P]$, $P \in \Pc$, or $U = (u_1, \dots, u_n)$, $n \ge 1$, $u_i \in
\Q$, and the following conditions hold:
\begin{itemize}\item no $u_i$ is an identity arrow, \item $u_i$
and $u_{i+1}$ do not belong to the same groupoid $\Hc$ or
$\Vc$.\end{itemize}

For instance, the horizontal identity $\id_{\Hc}P \in \Hc$, $P \in
\Pc$, is a path which is not reduced.

Every path $(u_1, u_2, \dots, u_n)$, with $n > 0$, can be
transformed into a reduced path by means of a finite number of
reductions, that is, operations of one of the following types:
\medbreak \begin{itemize}\item removing $u_i$ if $u_i$ is an
identity arrow and $n > 1$,

\medbreak
\item replacing $(u_1)$ by $[P]$ if $u_1 = \idh P$ or $\idv P$, $P \in
\Pc$,

\medbreak
\item replacing $(u_1, \dots, u_i, u_{i+1}, \dots, u_n)$ by
$(u_1, \dots, u_iu_{i+1}, \dots, u_n)$, if  $u_i$ and $u_{i+1}$
belong to the same groupoid $\Hc$ or $\Vc$.\end{itemize}

These operations generate an equivalence relation in $\path(\Q)$.
Following the lines of the proof of \cite[Theorem 5, Chapter
11]{higgins}, it is possible to see that in any equivalence class
there is a unique reduced path.

The set of all reduced paths on $\Q$ forms a groupoid under the
operation of concatenation followed by reduction. Compare also
with the analogous construction for groups \cite[p. 186]{s}.

Also, the set of all reduced paths on $\Q$ with this product is
isomorphic to $\Vc * \Hc$. Clearly $\Vc * \Hc$ contains both $\Vc$
and $\Hc$ as wide subgroupoids. In conclusion, any element $u$ of
$\Vc * \Hc$ has a unique standard form, namely:
\begin{itemize}
    \item $u \in \Pc$ (elements of length 0), or
    \item  $u=u_1u_2 \dots u_n$, where the $u_i$'s belong alternatively to different groupoids
    $\Vc$ or $\Hc$,  no $u_i$ is an identity (elements of length $n>0$).
\end{itemize}

In such case we shall say that $u_1$, respectively $u_n$, is the
\emph{first}, respectively the \emph{last}, letter of $u$.

\begin{lema}\label{lenghtadditive} Let $p= p_1 \dots p_N$, $q = q_1 \dots q_M$ be reduced paths in $\Vc * \Hc$
of lengths $N$ and $M$ respectively.

(i). If $p_N$ and $q_1$ belong to different groupoids $\Vc$ or
$\Hc$, then
 length $(pq) = N+M$.

(ii). If $p_N$ and $q_1$ belong to the same groupoid $\Vc$ or
$\Hc$, but $p_N\neq (q_1)^{-1}$ then
  length $(pq) = N+M - 1$.

(iii). If $p_N$ and $q_1$ belong to the same groupoid $\Vc$ or
$\Hc$, $p_N = (q_1)^{-1}$ but  $p_{N-1}\neq (q_2)^{-1}$ then
length $(pq) = N+M - 2$. \qed \end{lema}

\subsection{Diagonal groupoid of a slim double groupoid}\label{diagonalofslim}

\

Let $\T$ be a double groupoid. In the free product $\Vc * \Hc$, we
denote
$$
[A] := xgy^{-1}h^{-1}, \quad \text{if} \quad A =
\begin{matrix} \quad x \quad \\ h \,\, \boxe \,\, g \\
\quad y \quad \end{matrix} \in \B.
$$

Define the
'diagonal' groupoid $\D(\T)$ to be the quotient of the free
product $\Vc
* \Hc$ modulo the relations $[A]$, $A\in \B$.

\medbreak \emph{In the rest of this section we suppose that $\B$
is slim and satisfies the filling condition \ref{filling}}.

\begin{lema}\label{rel-norm} The subgroupoid $J$ generated by all relations $[A]$, $A\in \B$,
is a normal subgroup bundle of the free product $\Vc * \Hc$. \end{lema}

Hence, if $\B$ satisfies \ref{filling}, then $\D(\T) = (\Vc * \Hc)
/ J$.

\pf It is clear that $J$ is a subgroup bundle of $\Vc * \Hc$.
Let $A = \begin{matrix} \quad x \quad \\ h \,\, \boxe \,\, g \\
\quad y\quad
\end{matrix}$ be a box in $\B$.
We shall show that the expressions $z(xgy^{-1}h^{-1})z^{-1}$ and
$f(xgy^{-1}h^{-1})f^{-1}$ both belong to $J$, for all $z \in \Hc$,
$f \in \Vc$, such that $r(z) = l(x) = b(f)$. This implies
normality because $\Vc$ and $\Hc$ generate $\Vc * \Hc$. Let $z$,
$f$ as above. We have $r(z) = l(x) = t(h)$. Hence, since $\B$
satisfies \ref{filling}, we may pick a box $\begin{matrix} \quad z
\quad \\ r \,\, \boxe \,\, h \\ \quad s\quad \end{matrix}$ in
$\B$. Then the horizontal composition $$\begin{matrix} \quad z
\quad \\ r \,\, \boxe \,\, h \\ \quad s\quad \end{matrix} \,
\begin{matrix} \quad x \quad \\ h \,\, \boxe \,\, g \\ \quad
y\quad \end{matrix} = \begin{matrix} \quad zx \quad \\ r \,\,
\boxe \,\, g \\ \quad sy\quad \end{matrix}$$ is in $\B$.
Therefore, the expressions $X = zhs^{-1}r^{-1}$ and $Y =
zxgy^{-1}s^{-1}r^{-1}$ both belong to $J$. Hence so does the
product $YX^{-1} = z(xgy^{-1}h^{-1})z^{-1}$.
To show that $f(xgy^{-1}h^{-1})f^{-1}$ belongs to $J$, we argue as
before, now picking a box $B = \begin{matrix} \quad u \quad \\ f
\,\, \boxe \,\, v \\ \quad x\quad \end{matrix}$ in $\B$ and then
taking the vertical composition $\begin{matrix}B \vspace{-4pt} \\
A \end{matrix}$ in $\B$.  \epf

Composing the inclusions  $\Hc, \Vc \to \Vc * \Hc$ with the
canonical projection $\Vc * \Hc \to \D$ we get canonical groupoid
maps $i: \Hc \to \D$, $j: \Vc \to \D$.

\medbreak The diagram $(\D(\T), j, i)$ is  characterized by the
following universal property: for every diagram $(\G, j_0, i_0)$
over $\Vc$ and $\Hc$, such that $i_0(x)j_0(g) = j_0(h)i_0(y)$,
whenever the box $\begin{matrix} \quad x \quad \\ h \,\, \boxe \,\, g \\
\quad y\quad \end{matrix}$ is in $\B$, there is a unique morphism
of diagrams $f: \D(\T) \to \G$.

\medbreak It is clear that the assignment $\B \to (\D(\B), j, i)$
is functorial.

\medbreak The filling  condition on $\T$ corresponds to the
factorizability condition on $\D(\T)$, as we show next.

\begin{lema}\label{fill-fact}  $\D(\T) = j(\Vc) i(\Hc)$.
\end{lema}

In particular, if $\B$ is finite then $\D(\T)$ is a \emph{finite} groupoid.

\pf Every element in $\Vc * \Hc$  writes as a product $w_1 \dots
w_m$, where each $w_i$ is an element of $\Hc$ or $\Vc$. Hence
every element in $\D$ factorizes as $\overline{w_1} \dots
\overline{w_m}$, where $\overline{w_i}$ is the image of $w_i$
under the canonical projection, which coincides with $i(w_i)$ or
$j(w_i)$ according to whether $w_i$ is an element of $\Hc$ or
$\Vc$.

Therefore it is enough to see that every product $i(x)j(g)$, $x
\in \Hc, g \in \Vc$, belongs to $j(\Vc)i(\Hc)$. Indeed, this
implies that the elements in the factorization may be
appropriately reordered to get an element in $j(\Vc)i(\Hc)$. To
see this we use the assumption of condition \ref{filling} on $\T$:
since the product $i(x)j(g)$, $x \in \Hc, g \in \Vc$, is defined,
then there is a box $\begin{matrix} \quad x \quad \\ h \,\, \boxe
\,\, g \\ \quad y\quad \end{matrix}$ in $\B$. By construction of
$\D$, $i(x)j(g) = j(h)i(y)$. The lemma follows. \epf

\subsection{Main result}\label{mainresult}

\

We can now prove the main result of this section. Lemma
\ref{tecnico} encapsulates the most delicate part of the proof.

\begin{theorem}\label{class-slim} The assignments $\T \mapsto \D(\T)$ and $\D \mapsto \square(\D, j, i)$ determine
mutual category equivalences between
\begin{itemize}\item[(a)] The category of slim double groupoids $\begin{matrix} \B &\rightrightarrows &\Hc
\\\downdownarrows &&\downdownarrows \\ \Vc &\rightrightarrows &\Pc \end{matrix}$
 satisfying the filling condition \ref{filling}, with fixed $\Vc$ and $\Hc$, and

 \medbreak
\item[(b)] The category of $(\Vc, \Hc)$-factorizations of groupoids $\D$ on $\Pc$.
\end{itemize} \end{theorem}

\pf It remains to show that the assignments are mutually inverse.
Suppose first that $\D = j(\Vc) i(\Hc)$ is a factorization as in
(b). Let $\D'$ be the diagonal groupoid associated to the double
groupoid  $\square(\D, j, i)$; so that there are groupoid maps
$j': \Vc \to \D'$, $i': \Hc \to \D'$ such that $\D' = j'(\Vc)
i'(\Hc)$, in view of Lemma \ref{fill-fact}.

The universal property of $\D'$  implies the existence of a unique
groupoid map $f: \D' \to \D$ such that $fi' = i$ and $fj' = j$.
Moreover, $f$ is surjective because of the condition $\D = j(\Vc)
i(\Hc)$. To prove injectivity of $f$, let $P \in \Pc$ and let $z
\in \D'(P)$ such that $f(z) \in \Pc$. Write $z = j'(g)i'(x)$, $g
\in \Vc$, $x \in \Hc$, such that $b(g) = l(x)$. Applying $f$ to
this identity, we get $\id_P = f(z) = j(g)i(x)$. In particular
$t(g) = P = r(x)$, and $i(\id_P)j(\id_P) = j(g)i(x)$ in $\D$.

The definition of $\square(\D, j, i)$ implies that the box
$\begin{matrix} g \hspace{-0.5cm} &
\begin{tabular}{|p{0,1cm}||} \hhline{|=||} \\ \hline \end{tabular}
\\  & x\end{matrix}$ is in $\square(\D, j, i)$. Therefore, in view of the
defining relations in $\D'$, we have $z = j'(g)i'(x) = \id_P$.
This proves that $f$ is injective and thus an isomorphism.

\medbreak Let now $\T$ be a double  groupoid as in (a), $\D =
\D(\T)$ the associated diagonal groupoid with the canonical maps
$i: \Hc \to \D$, $j: \Vc \to \D$, and $\B' = \square(\D, j, i)$.
Let $\begin{matrix} \quad x \quad \\ h \,\, \boxe \,\, g \\ \quad
y\quad \end{matrix}$ be a box in $\B$. Then $i(x)j(g) = j(h)i(y)$
in $\D$ and this relation determines a unique box $\begin{matrix}
\quad x \quad \\ h \,\, \boxe \,\, g \\ \quad y\quad \end{matrix}$
in $\B'$, since $\B'$ is slim. This defines a map $F: \B \to \B'$
that, because of the slim condition on $\B$, turns out to be an
injective map of double groupoids.

We claim that $F$ is also surjective, hence an isomorphism. To
establish this claim, we shall need the presentation of the free
product $\Vc * \Hc$ given in Subsection \ref{presentation}.
Let $A = \begin{matrix} \quad x \quad \\ h \,\, \boxe \,\, g \\
\quad y \quad \end{matrix}$ be a box in $\B'$, which means that
\begin{equation}\label{duno}
i(x)j(g) = j(h)i(y)
\end{equation}
in $\D$.
 We shall prove that there is a box
$\begin{matrix} \quad x \quad \\ h \,\, \boxe \,\, g \\ \quad y
\quad \end{matrix}$ in $\B$. First, using the filling condition
\ref{filling} in $\B$, there
is a box $A_0 = \begin{matrix} \quad x \quad \\ h_0 \,\, \boxe \,\, g \\
\quad y_0 \quad \end{matrix} \in \B$. Then it is enough to show
that the box $E = \begin{matrix} f \hspace{-0.2cm} &
\begin{tabular}{|p{0,1cm}||} \hhline{|=||} \\ \hline \end{tabular}
\\  & z\end{matrix}$ is also in $\B$, where $f = h_0^{-1}h$ and $z = yy_0^{-1}$.
In fact, if  $E \in \B$, then $E \rightharpoondown
A_0 \in \B$ is the desired box.

Since  $A_0$ belongs to $\B$, $i(x)j(g) =
j(h_0)i(y_0)$ in $\D$; combined with \eqref{duno},
this gives $j(f)i(z) \in \Pc$. The
definition of $\D$ combined with Lemma \ref{rel-norm} implies that
the path $fz$ belongs to the normal group bundle $J$.
The proof of the Theorem will be achieved once the following Lemma is established.
\epf

\begin{lema}\label{tecnico} Let $f\in \Vc$ and  $z\in \Hc$ such that

\begin{itemize}
    \item $t(f) = r(z) =:P$ and $b(f) = l(z)$.
    \item There exist $A_1, \dots, A_n\in \B$, $\epsilon_1, \dots, \epsilon_n\in \{\pm 1\}$
    such that
    \begin{equation}\label{exp-red}
    fz = [A_1]^{\epsilon_1} \dots [A_n]^{\epsilon_n}.\end{equation}

Then there exists $E\in \Ec$ such that $E = \begin{matrix} f
\hspace{-0.2cm} &
\begin{tabular}{|p{0,1cm}||} \hhline{|=||} \\ \hline \end{tabular}
\\  & z\end{matrix}$.
\end{itemize}
\end{lema}

\pf We proceed by induction on $n$. Note
that the case $n=0$ means that $fz = \id_P$ in $\Vc * \Hc$, so
that $f = \id^{\Vc}_P$, $z = \id^{\Hc}_P$, since the word $fz$ is
not reduced; thus $E = \Theta_P$ is the desired element in $\Pc$.

\medbreak {\it Proof when $n=1$.} We need to work out this case by
technical reasons, see Sublemma \ref{zutecnico} below. We shall
assume that $f\notin \Pc$, $z\notin \Pc$, so that the path $fz$ is
reduced; if either $f\in \Pc$ or $z\in \Pc$ the arguments are
similar. We have $fz = [A_1]^{\epsilon_1}$. If $ \epsilon_1 = 1$
this says that $fz = xgy^{-1}h^{-1}$ where we omit the subscript 1
for simplicity. Hence, the right-hand side is not reduced. Since
the letters $x$, $g$, $y$, $h$ belong alternatively to different
groupoids, at least one of them should be in $\Pc$. If $x = \idh
P$, then  $fz = gy^{-1}h^{-1}$ hence necessarily $h = \idv P$,
$f=g$, $z=y^{-1}$ and $A_1 =
\begin{matrix}
\begin{tabular}{||p{0,1cm}|} \hhline{||=|} \\ \hline \end{tabular}
\,f \\  z^{-1}  \end{matrix}$. Thus $E = (A_1)^h$ does the job. If
$x \notin \Pc$ then it should be cancelled because the left-hand
side begins by $f\in \Vc$; thus $g\in \Pc$, $y=x^{-1}$ and $fz =
h$, a contradiction.

\medbreak Assume now that $ \epsilon_1 = -1$, that is, $fz =
hyg^{-1}x^{-1}$. Again, the right-hand side is not reduced and one
of the letters should be in $\Pc$. If $h = \idv P$ then $fz =
yg^{-1}x^{-1}$ hence necessarily $y = \idh P$, $f=g^{-1}$,
$z=x^{-1}$ and  $E = (A_1)^{-1}$ does the job. If $h \notin \Pc$
then either
\begin{itemize}
    \item [(a)] $y\in \Pc$, $f=hg^{-1}$, $z = x^{-1}$ and
    $E = \left\{\begin{matrix} \id h \\ (A_1)^{-1} \end{matrix}\right\}$ does the job; or
    \item [(b)] $y\notin \Pc$, $g\in \Pc$, $f = h$, $z = yx^{-1}$ and
    $E =   A_1\, \id x^{-1}$ does the job.
\end{itemize}

\medbreak {\it Assume now that the claim is true for $n-1$.}
Assume that \eqref{exp-red} holds for $n$. Our aim is to reduce
the right-hand side of \eqref{exp-red} to $n-1$ factors, or to
achieve a  contradiction by comparison of the lengths in both
sides.

Before we begin to analyze contiguous brackets, where `bracket'
means an element of the form $[A_i]^{\pm 1}$, let us set up some
preliminaries. Let us say that $A_i\in \B$ is of \emph{class
$\ell$} if exactly $\ell$ of its sides are not in $\Pc$. We shall
assign a `type' to any bracket $[A_i]^{\pm 1}$.

\begin{itemize}
    \item If $\ell =0$, all the sides of $A_i$ are in $\Pc$; hence
    $[A_i]^{\epsilon_i}$ can be extirpated from the right-hand
    side of \eqref{exp-red} and we are done by the inductive
    hypothesis. Hence, we can assume that $\ell >0$ for all $i$.
\end{itemize}

\medbreak\begin{itemize}
    \item If $\ell = 1$, we distinguish two types.
\end{itemize}

($\alpha$): $A_i$ has a non-trivial vertical side. We can assume
that $A_i = h_i
\begin{tabular}{|p{0,1cm}||} \hhline{||=||} \\ \hhline{||=||}
\end{tabular}$. For, if $A_i =
\begin{tabular}{||p{0,1cm}|} \hhline{||=||} \\ \hhline{||=||}
\end{tabular}\; g_i$ then $[A_i] = [\widetilde A_i]$ where
$\widetilde A_i = \left\{\begin{matrix} \id g_i^{-1} \\ A_i
\end{matrix}\right\} = g_i^{-1}
\begin{tabular}{|p{0,1cm}||} \hhline{||=||} \\ \hhline{||=||}
\end{tabular}$.

Moreover if $[A_i]$ is of type $(\alpha)$ then $[A_i]^{-1} =
[A_i^{v}]$, again of type $(\alpha)$.

\medbreak ($\beta$): $A_i$ has a non-trivial horizontal side. We
can assume that $A_i = \begin{matrix} x_i\\
\begin{tabular}{||p{0,1cm}||} \hline \\ \hhline{||=||}
\end{tabular}\\ \quad\end{matrix}$.

\noindent Moreover if $[A_i]$ is of type $(\beta)$ then
$[A_i]^{-1} = [A_i^{h}]$, again of type $(\beta)$. In particular,
if $\ell = 1$, then we can assume $\epsilon_i = 1$.

\medbreak
\begin{itemize}
    \item If $\ell = 2$, we can assume that the two non-trivial sides
    live in different groupoids. For, if $A_i = h_i
\begin{tabular}{|p{0,1cm}|} \hhline{||=||} \\ \hhline{||=||}
\end{tabular}\, g_i$ then $[A_i] = [\widetilde A_i]$ where
$\widetilde A_i = \left\{\begin{matrix} \id g_i^{-1} \\ A_i
\end{matrix}\right\} = g_i^{-1}h_i
\begin{tabular}{|p{0,1cm}||} \hhline{||=||} \\ \hhline{||=||}
\end{tabular}$, whose bracket is of type $(\alpha)$. Similarly, if $A_i$ has two
non-trivial horizontal sides then $[A_i]^{\epsilon_i}$ can be
replaced by a  bracket of type $(\beta)$.
\end{itemize}

\medbreak We distinguish two types.

($\gamma$): $[A_i]\in \Hc\Vc$. We can assume that $A_i =
\begin{matrix} \quad\; x_i \\h_i
\begin{tabular}{|p{0,1cm}||} \hline \\ \hhline{||=||}
\end{tabular}\\ \quad \end{matrix}$. For, if $A_i =\begin{matrix}  x_i\quad \\
\begin{tabular}{||p{0,1cm}|} \hline \\ \hhline{||=||}
\end{tabular}\; g_i\\ \quad \end{matrix}$ then $[A_i] = [\widetilde A_i]$ where
$\widetilde A_i = \left\{\begin{matrix}   A_i \\ \id g_i^{-1}
\end{matrix}\right\} =
\begin{matrix} \qquad x_i \\g_i^{-1}
\begin{tabular}{|p{0,1cm}||} \hline \\ \hhline{||=||}
\end{tabular}\\ \quad \end{matrix}$. Similarly, if $A_i =\begin{matrix}
h_i\begin{tabular}{|p{0,1cm}||} \hhline{||=||} \\ \hline
\end{tabular}\\\quad  y_i\end{matrix}$ then $[A_i] = [\widetilde A_i]$ where
$\widetilde A_i =    A_i  \id y_i^{-1} =
\begin{matrix} \qquad y_i^{-1} \\h_i
\begin{tabular}{|p{0,1cm}||} \hline \\ \hhline{||=||}
\end{tabular}\\ \quad \end{matrix}$.

\medbreak ($\delta$): $[A_i]\in \Vc\Hc$. Then $A_i =\begin{matrix}
\begin{tabular}{||p{0,1cm}|} \hhline{||=||} \\ \hline
\end{tabular} \; g_i\\y_i \quad  \end{matrix}$.

\medbreak \noindent Moreover if $[A_i]$ is of type $(\gamma)$ then
$[A_i]^{-1} = [B]$, where $B = (A_i \id t(A_i)^{-1})^h$ is of type
$(\delta)$. In particular, if $\ell = 2$ then we can assume
$\epsilon_i = 1$.

\medbreak \begin{itemize}
    \item If $\ell = 3$, we can assume that the identity sides
    are either $x_i$ or $h_i$;
    otherwise we replace the box by one with $\ell =2$.
    There are two types.
\end{itemize}
\begin{equation*} (\zeta): \quad A_i =\begin{matrix} x_i \quad \\
\begin{tabular}{||p{0,1cm}|} \hline \\ \hline
\end{tabular} \; g_i\\y_i \quad  \end{matrix}, \qquad
(\eta): \quad A_i =\begin{matrix} h_i
\,\begin{tabular}{|p{0,1cm}|} \hhline{||=||} \\ \hline
\end{tabular} \; g_i\\y_i   \end{matrix}.
\end{equation*}

\medbreak Moreover if $[A_i]$ is of type $(\zeta)$ then
$[A_i]^{-1} = [A_i^{v}]$ is again of type $(\zeta)$; if $[A_i]$ is
of type $(\eta)$ then $[A_i]^{-1} = [A_i^{h}]$ is again of type
$(\eta)$. In particular, if $\ell = 3$ then we can assume
$\epsilon_i = 1$.

\begin{itemize}
    \item If $\ell = 4$, we distinguish two types.
\end{itemize}
\begin{equation*} (\theta): \quad\epsilon = 1, \qquad(\kappa): \quad\epsilon = -1.
\end{equation*}

\medbreak \emph{In conclusion  the right-hand side of
\eqref{exp-red} is a product of $n$ brackets of types $(\alpha)$,
\dots, $(\kappa)$ with the exponent $\epsilon_i = 1$ except for
the type $(\kappa)$.}


\medbreak Next, we can assume several restrictions on the
contiguity of these brackets, as summarized in the following
statement. We prove below that, whenever these restrictions do not
hold, then we can replace the pair of contiguous brackets by a
single bracket and hence apply the inductive hypothesis.

\medbreak
\begin{zulma}\label{restrictions} Let $1\le i \leq n$.

\medbreak \noindent{\rm (1).} If $[A_{i}]$ is of type $(\alpha)$,
then $i>1$ and $[A_{i-1}]^{\epsilon_{i-1}}$ is of type $(\kappa)$.
Also, if $i<n$ then $[A_{i+1}]^{\epsilon_{i+1}}$ is of type
$(\theta)$.

\medbreak \noindent{\rm (2).} If $[A_{i}]$ is of type $(\beta)$,
then  $i<n$ and $[A_{i+1}]^{\epsilon_{i+1}}$ is of type
$(\kappa)$.
    Also,  if $i>1$ then $[A_{i-1}]^{\epsilon_{i-1}}$ is of type $(\theta)$.

\medbreak \noindent{\rm (3).} If $[A_{i}]$ is of type $(\gamma)$,
then:
\begin{description}
    \item[(a)] If $i>1$ then $[A_{i-1}]^{\epsilon_{i-1}}$ is of type  $(\eta)$, $(\theta)$,
    or $(\kappa)$ with $x_{i-1} \neq x_i$.
    \item[(b)] If $i<n$ then $[A_{i+1}]^{\epsilon_{i+1}}$ is of type  $(\zeta)$, $(\theta)$,
    or $(\kappa)$ with $h_{i+1} \neq h_i$.
\end{description}

\medbreak \noindent{\rm (4).} If $[A_{i}]$ is of type $(\delta)$,
then:
\begin{description}
    \item[(a)] If $i>1$ then $[A_{i-1}]^{\epsilon_{i-1}}$ is of type $(\zeta)$, $(\theta)$,
    $(\kappa)$, or $(\theta)$ with $h_{i-1} \neq g_i$.
    \item[(b)] If $i<n$ then $[A_{i+1}]^{\epsilon_{i+1}}$ is of type $(\eta)$, $(\kappa)$,
    or $(\theta)$ with $x_{i+1} \neq y_i$.
\end{description}

\medbreak \noindent{\rm (5).} If $[A_{i}]$ is of type $(\zeta)$,
then:
\begin{description}
    \item[(a)] If $i>1$ then $[A_{i-1}]^{\epsilon_{i-1}}$ is of type $(\gamma)$,
    $(\zeta)$, $(\eta)$, $(\theta)$ or $(\kappa)$.
Also,  \begin{enumerate}
    \item[(1)] if $[A_{i-1}]$ is of type $(\zeta)$ then $y_{i-1} \neq x_i$;
    \item[(2)] if $[A_{i-1}]^{\epsilon_{i-1}}$ is of type $(\kappa)$ then $x_{i-1} \neq x_i$.
\end{enumerate}

    \item[(b)] If $i<n$ then $[A_{i+1}]^{\epsilon_{i+1}}$ is of type
    $(\delta)$, $(\zeta)$,  $(\eta)$, $(\theta)$ or $(\kappa)$.
Also,
\begin{enumerate}
    \item[(1)] if $[A_{i+1}]$ is of type $(\zeta)$ then $y_{i} \neq x_{i+1}$;
    \item[(2)] if $[A_{i+1}]^{\epsilon_{i+1}}$ is of type $(\theta)$ then $y_{i} \neq x_{i+1}$.
\end{enumerate}
\end{description}

\medbreak \noindent{\rm (6).} If $[A_{i}]$ is of type $(\eta)$,
then:
\begin{description}
    \item[(a)] If $i>1$ then $[A_{i-1}]^{\epsilon_{i-1}}$ is of type  $(\delta)$, $(\zeta)$,  $(\eta)$, $(\theta)$ or $(\kappa)$.
Also,
\begin{enumerate}
    \item[(1)] if $[A_{i-1}]$ is of type $(\eta)$ then $h_{i-1} \neq g_i$;
    \item[(2)] if $[A_{i-1}]^{\epsilon_{i-1}}$ is of type $(\theta)$ then $h_{i-1} \neq g_i$.
\end{enumerate}

    \item[(b)] If $i<n$ then $[A_{i+1}]^{\epsilon_{i+1}}$ is of type $(\gamma)$, $(\zeta)$, $(\eta)$, $(\theta)$ or $(\kappa)$.
Also,
\begin{enumerate}
    \item[(1)] if $[A_{i+1}]$ is of type $(\eta)$ then $g_{i+1} \neq h_i$;
    \item[(2)] if $[A_{i+1}]^{\epsilon_{i+1}}$ is of type $(\kappa)$ then $h_{i+1} \neq h_i$.
\end{enumerate}
\end{description}

\medbreak \noindent{\rm (7).} If $[A_i]^{\epsilon_i}$ and
$[A_{i+1}]^{\epsilon_{i+1}}$ are both of types $(\theta)$ or
$(\kappa)$ for some $i$, $1\le i \le
    n-1$, then either

\begin{itemize}
    \item     $\epsilon_i =\epsilon_{i+1}$; or
    \item $\epsilon_i = 1 = -\epsilon_{i+1}$,
    and $h_i \neq h_{i+1}$; or
    \item $\epsilon_i = 1 = -\epsilon_{i+1}$,
    $h_i = h_{i+1}$ and $y_i \neq y_{i+1}$; or
    \item $\epsilon_i = -1 = -\epsilon_{i+1}$ and $x_i \neq x_{i+1}$; or
    \item $\epsilon_i = -1 = -\epsilon_{i+1}$, $x_i = x_{i+1}$
    and $g_i \neq g_{i+1}$.
\end{itemize}
    \end{zulma}

\medbreak We deal with (1). If $[A_1]$ is of type $(\alpha)$ then
$ h_1fz = [A_2]^{\epsilon_2} \dots [A_n]^{\epsilon_n}$; by the
inductive hypothesis there exists $E_1\in \Ec$ such that $E_1 =
\begin{matrix} h_1f \hspace{-0.2cm} &
\begin{tabular}{|p{0,1cm}||} \hhline{|=||} \\ \hline \end{tabular}
\\  & z\end{matrix}$. Then $E = \left\{\begin{matrix} (A_1)^v\\ E_1 \end{matrix}\right\}$
does the job. Thus we do not consider this possibility. Assume
that $i>1$ and $[A_{i-1}]^{\epsilon_{i-1}}$ is not of type
$(\kappa)$. We show that $[A_{i-1}][A_i] = [B]$; hence the
right-hand side of \eqref{exp-red} has really $n-1$ factors and
the existence of $E$ follows from the inductive hypothesis. Thus
we do not consider this possibility. Explicitly, $B\in \B$ is
given as follows.

\begin{itemize}
    \item If $[A_{i-1}]$ is of type $(\alpha)$ or $(\eta)$ then $B = \left\{\begin{matrix} A_i\\ A_{i-1} \end{matrix}\right\}$.
    \item If $[A_{i-1}]$ is of type $(\beta)$ then $B =  \left\{\begin{matrix} A_{i-1}\\ A_i \end{matrix}\right\}$.
    \item If $[A_{i-1}]$ is of type $(\gamma)$ then $B =
    \left\{\begin{matrix} A_i\\ \id h_{i-1}^{-1}\\ A_{i-1} \end{matrix}\right\}$.
    \item If $[A_{i-1}]$ is of type $(\delta)$ or $(\zeta)$ then $B = A_i\, A_{i-1}$.
    \item If $[A_{i-1}]$ is of type $(\theta)$ then
    $B = \left\{\begin{matrix} A_i\\ A_{i-1}C \end{matrix}\right\}C^h$, where $C\in \B$ is a box filling
    the configuration $\begin{matrix} \qquad x_i^{-1}\\ g_i\begin{tabular}{|p{0,3cm}}  \hline \\\end{tabular}
    \end{matrix}$.
            \end{itemize}

\medbreak Now, if $i<n$ then write $[A_{i+1}]^{\epsilon_{i+1}} =
[\widetilde A_{i+1}]^{-\epsilon_{i+1}}$. Then
$[A_{i}][A_{i+1}]^{\epsilon_{i+1}} = \left([\widetilde
A_{i+1}]^{-\epsilon_{i+1}}[A_{i}^v]\right)^{-1}$. We discard this
possibility by the previous discussion unless $[\widetilde
A_{i+1}]^{-\epsilon_{i+1}}$ is of type $(\kappa)$, that is, unless
$[A_{i+1}]^{\epsilon_{i+1}}$ is of type $(\theta)$.

\medbreak We deal with  (2). If $[A_n]$ is of type $(\beta)$ then
$ fzx_n^{-1} = [A_1]^{\epsilon_1} \dots
[A_{n-1}]^{\epsilon_{n-1}}$; by the inductive hypothesis there
exists $E_1\in \Ec$ such that $E_1 =
\begin{matrix} f \hspace{-0.2cm} &
\begin{tabular}{|p{0,1cm}||} \hhline{|=||} \\ \hline \end{tabular}
\\  & zx_n^{-1}\end{matrix}$. Then $E = E_1  (A_n)^v$
does the job. Thus we do not consider this possibility. Assume
that $i<n$ and $[A_{i+1}]^{\epsilon_{i+1}}$ is not of type
$(\kappa)$. We show that $[A_i][A_{i+1}] = [B]$; hence the
right-hand side of \eqref{exp-red} has really $n-1$ factors and
the existence of $E$ follows from the inductive hypothesis. Thus
we do not consider this possibility. Explicitly, $B\in \B$ is
given as follows.

\begin{itemize}
    \item If $[A_{i+1}]$ is of type $(\alpha)$ or $(\eta)$ then $B = \left\{\begin{matrix} A_i\\ A_{i+1} \end{matrix}\right\}$.
    \item If $[A_{i+1}]$ is of type $(\beta)$ or $(\delta)$ or $(\zeta)$ then $B =    A_i\, A_{i+1}$.
    \item If $[A_{i+1}]$ is of type $(\gamma)$ then $B =
A_i\left\{\begin{matrix}   A_{i+1}\\\id h_{i+1}^{-1}
\end{matrix}\right\}$.
    \item If $[A_{i+1}]$ is of type $(\theta)$ then
    $B = \left\{\begin{matrix} A_i\left\{\begin{matrix} A_{i+1}\\ C \end{matrix}\right\}\\ C^v
    \end{matrix}\right\}$,
    where $C\in \B$ is a box filling
    the configuration $\begin{matrix} \qquad y_{i+1}\\ h_{i+1}^{-1}\begin{tabular}{|p{0,3cm}}  \hline \\\end{tabular}
    \end{matrix}$.
            \end{itemize}

If $i>1$, we may write $[A_{i-1}]^{\epsilon_{i-1}} = [\widetilde
A_{i-1}]^{-\epsilon_{i-1}}$, hence
$[A_{i-1}]^{\epsilon_{i-1}}[A_{i}]= \left([A_{i}^h][\widetilde
A_{i-1}]^{-\epsilon_{i-1}}\right)^{-1}$. We discard this, unless
$[A_{i-1}]^{\epsilon_{i-1}}$ is of type $(\theta)$.

\medbreak We deal with  (3a). We may assume that
$[A_{i-1}]^{\epsilon_{i-1}}$ is neither of type $(\alpha)$ nor
$(\beta)$ by (1) and (2). If $[A_{i-1}]^{\epsilon_{i-1}}$ is of
type $(\gamma)$, $(\delta)$, $(\zeta)$, or $(\kappa)$ with
$x_{i-1} = x_i$, then $[A_{i-1}]^{\epsilon_{i-1}}[A_i] = [B]$;
hence the right-hand side of \eqref{exp-red} has $n-1$ factors
and, by the inductive hypothesis, we do not consider this
possibility. Here $B\in \B$ is:

\begin{itemize}
    \item If $[A_{i-1}]$ is of type $(\gamma)$ then $B =
\left\{\begin{matrix}   A_{i-1}\\\id h_{i-1}^{-1}
\\ \left\{A_{i}\id x_{i}^{-1}\right\}\end{matrix}\right\}$.
    \item If $[A_{i-1}]$ is of type $(\delta)$ or $(\zeta)$ then $B =
A_{i}\, \id x_{i}^{-1}\, A_{i-1}$.
\item If $[A_{i-1}]^{-1}$ is of type $(\kappa)$
with $x_{i-1} = x_i$, then $B =
\left\{\begin{matrix}   A_{i} \id x_{1}^{-1}
\\ A_{i-1}\end{matrix}\right\}$.
            \end{itemize}

\medbreak We deal with  (3b). We may assume that
$[A_{i+1}]^{\epsilon_{i+1}}$ is neither of type $(\alpha)$ nor
$(\beta)$ nor $(\gamma)$  by (1), (2) and (3a). If
$[A_{i+1}]^{\epsilon_{i+1}}$ is of type $(\delta)$, $(\eta)$, or
$(\kappa)$ with $h_i = h_{i+1}$ then
$[A_i][A_{i+1}]^{\epsilon_{i+1}} = [B]$. By the inductive
hypothesis, we do not consider this possibility. Here $B\in \B$
is:

\begin{itemize}
    \item If $[A_{i+1}]$ is of type $(\delta)$ or $(\eta)$ then $B = \left\{\begin{matrix}
A_{i}\\\id h_{i}^{-1} \\ A_{i+1}\end{matrix}\right\}$.
    \item If $[A_{i+1}]^{-1}$ is of type $(\kappa)$
with $h_i = h_{i+1}$  then $B =
\left\{\begin{matrix}   A_{i} \\ \id h_{i}^{-1}
\end{matrix}\right\}A_{i+1}^v$.
            \end{itemize}

\medbreak Now (4a) follows from (3b), and (4b) follows from (3a),
by the inversion argument as for the second parts of (1) and (2).

\medbreak We deal with  (5a). If $[A_{i-1}]$ is of type $(\zeta)$
and $y_{i-1} = x_i$, then $[A_{i-1}][A_i] = [B]$ where $B =
\left\{\begin{matrix} A_{i-1}\\A_{i}\end{matrix}\right\}$. If
$[A_{i-1}]^{-1}$ is of type $(\kappa)$ and $x_{i-1} = x_i$, then
$[A_{i-1}]^{-1}[A_i] = [B]^{-1}$ where $B = \left\{\begin{matrix}
A_{i}^v\\A_{i-1}\end{matrix}\right\}$. By the inductive
hypothesis, we do not consider this possibility. Now (5b) follows
from (5a) by the inversion argument.

\medbreak We deal with  (6a). If $[A_{i-1}]$ is of type $(\eta)$
and $h_{i-1} = g_i$, then $[A_{i-1}][A_i] = [B]$ where $B =
 A_{i}A_{i-1}$. If
$[A_{i-1}]$ is of type $(\theta)$ and $h_{i-1} = g_i$, then
$[A_{i-1}][A_i] = [B]$ where $B = A_{i} A_{i-1}$. By the inductive
hypothesis, we do not consider this possibility. Now (6b) follows
from (6a) by the inversion argument.

\medbreak We deal with  (7). If  $\epsilon_i = 1 =
-\epsilon_{i+1}$, $h_i = h_{i+1}$
    and $y_i = y_{i+1}$, then $[A_{i}][A_{i+1}]^{-1} = [B]$ where $B = \left\{\begin{matrix}
A_i \\ A_{i+1}^v\end{matrix}\right\}$.

If $\epsilon_i = -1 = -\epsilon_{i+1}$, $x_i = x_{i+1}$
    and $g_i = g_{i+1}$, then
    $[A_{i}]^{-1}[A_{i+1}] = [B]^{-1}$ where $B = A_iA_{i+1}^h$.

To conclude the proof of the lemma, and {\it a fortiori} of the
theorem, we need to establish the following fact.

\begin{zulema}\label{zutecnico}
Let $n\ge 2$. Consider an element $P_n \in \Vc*\Hc$, such that
    \begin{equation}\label{exp-redbis}
    P_n = [A_1]^{\epsilon_1} \dots [A_n]^{\epsilon_n},\end{equation}
where $A_1, \dots, A_n\in \B$, $\epsilon_1, \dots, \epsilon_n\in
\{\pm 1\}$; the brackets $[A_i]^{\epsilon_i}$ are of type
$(\alpha)$, \dots, $(\kappa)$, $1\le i \le n$; and contiguous
brackets satisfy Restrictions \ref{restrictions}. Then
\begin{flalign*}
 &\text{(i) the last letter of }    P_n \text{ is the
last letter of } [A_n]^{\epsilon_n},&
\\ &\text{(ii) if } [A_n]^{\epsilon_n} \text{has type $(\theta)$
or $(\kappa)$ then the last two letters of }    P_n
\intertext{equal the last two of  $[A_n]^{\epsilon_n}$,  and}
&\text{(iii) the length of }    P_n > 2.&
\end{flalign*}

\end{zulema}
\pf If $n = 2$ the claim follows by inspection of the possible
cases; in fact, the length of $P_2$ turns out to be $\ge 5$. Now
assume that the claim is true for $n\ge 2$. Consider $P_{n+1}$ as
in \eqref{exp-redbis} and write $P_{n+1} = P_n
[A_{n+1}]^{\epsilon_{n+1}}$. By hypothesis, the length of $P_2$ is
$>2$ and its last letter is that of $[A_n]^{\epsilon_n}$. Now the
pair $[A_n]^{\epsilon_n}[A_{n+1}]^{\epsilon_{n+1}}$ satisfies
Restrictions \ref{restrictions}. The sublemma now follows from
Lemma \ref{lenghtadditive}.
 \epf

We can now finish the proof of the Lemma. If the right-hand side
of \eqref{exp-red} satisfies Restrictions \ref{restrictions}, then
the sublemma gives a contradiction, since the left-hand side has
length $\le 2$. Thus at least one bracket can be eliminated in the
right-hand side and the Lemma follows by induction. \epf

\section{Fusion double groupoids}\label{fusion}

Let  $\T$ be a finite double groupoid satisfying the filling
condition \ref{filling}. The following definition is motivated by
\cite[Proposition 3.16]{AN2}.

\begin{definition}\label{defi:fusion}
We say that $\T$ is \emph{fusion} if and only if the following
hold:
\begin{enumerate}
\item[(F1)]\label{f1} The vertical groupoid $\Vc\rightrightarrows \Pc$ is connected.

\medbreak
\item[(F2)]\label{f2} For any $x\in \Hc$, there exists at most one $E\in \Ec$
such that $b(E) = x$.
\end{enumerate}
\end{definition}

Observe that condition (F2) in the definition is equivalent to
injectivity of the morphism $b: \Ec \to \Hc$ of groupoids over
$\Pc$.

\medbreak Let $\Ec$ be the core groupoid of $\B$ and let $\mathfrak E$ be
the core groupoid of its frame $\F$. Then there is an exact
sequence of groupoids $$\begin{CD}1 @>>> \Kc @>>> \Ec @>>>
\mathfrak E @>>> 1 \end{CD}.$$ In particular, for any $P\in \Pc$,
we have
\begin{equation}\label{core-slim}
\vert \Ec(P) \vert = \vert \Kc(P) \vert \vert \mathfrak E(P)
\vert.
\end{equation}

\begin{prop} Suppose that $\T$ is fusion. Then $\T$ is slim.
\end{prop}

\pf It follows immediately from condition (F2) in Definition
\ref{defi:fusion}. \epf

Then, in view of Theorem \ref{class-slim}, a fusion double
groupoid $\T$ is determined by a $(\Vc, \Hc)$-factorization of its
diagonal groupoid $\D = \D(\T)$.

\begin{prop} The following are equivalent:

\begin{itemize} \item[(i)] $\T$ is fusion. \item[(ii)] $\Vc \rightrightarrows \Pc$ is connected and $j: \Vc \to \D(\T)$ is injective.
\end{itemize}
In particular, if $\T$ is fusion, then its diagonal groupoid $\D=
\D(\T)$ is connected. \end{prop}

\pf It is enough to see that injectivity of the map $b: \Ec \to
\Hc$ is equivalent to injectivity of the map $j$. Suppose first
that $b: \Ec \to \Hc$ is injective. Let $g \in \Vc$ such that
$j(g) \in \Pc$. By Theorem \ref{class-slim}, $\T \simeq \T(\D, j,
i)$. Then there is a box
$g \, \begin{tabular}{|p{0,1cm}||} \hhline{|=|} \\
\hhline{|=|}\end{tabular} \in \B$. Since this box belongs to the
core groupoid $\Ec$, the injectivity of $b$ implies that $g \in
\Pc$. Hence $j$ is injective.

\medbreak Conversely, suppose that $j$ is injective. Let $E \in
\Ec$ such that $b(E) \in \Pc$, so that
$E =  g\, \begin{tabular}{|p{0,1cm}||} \hhline{|=|} \\
\hhline{|=|}\end{tabular}$, for some $g \in \Vc$. The construction
of $\D$ implies that $j(g) \in \Pc$ and therefore $g \in \Pc$. The
slim condition on $\T$ now guarantees that $E \in \Pc$. Thus $b$
is injective. \epf

\appendix\section{Corner functions}

\

\subsection{Orbits of groupoid actions}\label{orbits}

\ We first discuss an extension of the basic counting formula for
orbits in group theory  to the setting of groupoids. Let $s, e: \G
\rightrightarrows \Pc$ be a groupoid. Recall the equivalence
relation induced by $\G$ on $\Pc$: $P\sim Q$ if and only if $\G(P,
Q) \neq \emptyset$. We denote by $\widetilde Q$ the equivalence
class of $Q$. We set
$$
\GQ = \{g\in \G: e(g) = Q\} = \coprod_{P\in \Pc} \G(P, Q) =
\coprod_{P\in \widetilde Q} \G(P, Q),
$$
the set of arrows with target $Q$. It is clear from the above that
$$|\GQ| = |\widetilde Q| \times \vert \G(Q)\vert.$$ We shall
consider the fiber bundle $s: \GQ \to \Pc$.

\medbreak Let $K$ be a subgroup of $\G(Q)$. We define the quotient
$\GQ / K := \GQ / \equiv_K$, where $\equiv_K$ is the equivalence
relation in $\GQ$ given by
$$g \equiv_K h \iff g^{-1}h \in K.$$
Clearly, the source map descends to the quotient and we can
consider the fiber bundle $s: \GQ/K \to \Pc$. Hence, if $\G$ is
finite, then
$$|\GQ/K| = \dfrac{|\widetilde Q| \times \vert \G(Q)\vert}{\vert
K\vert}.$$ Clearly, $\G$ acts on $s: \GQ \to \Pc$ by left
multiplication.

\medbreak Assume that $\G$ acts on $p: \Ee \to \Pc$. The groupoid
$\G$ still acts on the orbit $\Oc_x$. Then there is an isomorphism
of $\G$-fiber bundles $\varphi: \Gpx / \G^x \to \Oc_x$ induced by
$g\mapsto g\fde x$. In particular,  if $\G$ is finite, then
\begin{equation}\label{conteo}
|\Oc_x| = \dfrac{|\widetilde{p(x)}| \times \vert
\G(p(x))\vert}{\vert \G^x\vert}.
\end{equation}

\subsection{Corner functions}\label{corners}

Let $\B$ be a finite double groupoid. We apply the counting
argument above to give alternative proofs of some properties of
the 'corner' functions defined in \cite{AN2}. There are four
corner functions but it is enough to consider one. Recall the sets
$\Ur(x,g)$ and $\Ur(B)$ with prescribed upper and right sides, as
defined in \eqref{ur}. The `upper-right' corner functions are
$\urcorner: \Hc \Times{r}{t}\Vc \to \mathbb N \cup \{ 0 \}$ and
$\urcorner: \B \to \mathbb N$ defined by
$$\urcorner (x,g) = |\Ur(x,g)|\quad\text{ and }\quad\urcorner (B)
= |\Ur(B)|.$$ The other corner functions are defined similarly. In
this parlance, the \emph{filling condition} on $\T$ is just
\begin{equation}\label{filling-app}
\urcorner (x,g) >0 \quad \text{for any }(x,g)\in \Hc
\Times{r}{t}\Vc.
\end{equation}

 We now interpret the corner functions in terms of orbits
of an action of the core groupoid.  Let $\gamma: \B \to \Pc$ be
the `left-bottom' vertex, $\gamma(B) = lb(B)$. Let $B\in \B$ and
$Q = \gamma(B) = bl(B)$. Consider the relation $\sim$ on $\Pc$
induced by $\Ec$. Recall that $\theta(Q)$ is the common value
\begin{equation*} \lrcorner(\id_{\Vc} Q, \id_{\Hc} Q)
= \urcorner(\id_{\Vc} Q, \id_{\Hc} Q) = \llcorner(\id_{\Vc} Q,
\id_{\Hc} Q) = \ulcorner(\id_{\Vc} Q, \id_{\Hc} Q).
\end{equation*}
The proposition \ref{prop-core}, together with formula
\eqref{conteo}, implies the following formula for the corner
function:
\begin{equation}\label{formulacorner} \urcorner (B) =
|\widetilde{Q}| \times \vert \Ec(Q)\vert.
\end{equation}
Applied to $B = \Theta_Q$, the formula implies that $$\theta(Q) =
\urcorner(\Theta_Q) = |\widetilde{Q}| \times \vert \Ec(Q)\vert =
\urcorner (B).$$ Hence $\theta(Q)$ is also given by
\eqref{formulacorner}. That is, the corner functions on a box
depend only on the vertex 'opposite' to the corner of that box.
Formula \eqref{formulacorner} provides easy alternative proofs of
the following facts-- see \cite{AN2}:

\medbreak
\begin{enumerate}
\item[(a)] Let $P, Q \in \Pc$. If $P\sim Q$, then $\theta(P) = \theta(Q)$.

\medbreak
\item[(b)]
Let $L, M, N \in \B$. Suppose that
\begin{tiny}\begin{tabular}{p{0,25cm}|p{0,3cm}} $L$ &
$M$ \\ \hline $N$ & \quad \end{tabular}\end{tiny}. Then
\begin{align*}  \quad \ulcorner(L) = \urcorner(M), \quad \llcorner(L) = \lrcorner(M), \quad \ulcorner(L) =
\llcorner(N), \quad \urcorner(L) = \lrcorner(N).\end{align*}

\medbreak
\item[(c)] Let $X, Y, Z \in \B$ such that
\begin{tiny}\begin{tabular}{p{0,25cm}|p{0,25cm}} $X$ & $Y$ \\
\hline $Z$ & \quad \end{tabular}\end{tiny}. Then
\begin{align*}
\urcorner(XY) = \urcorner(X), \quad \urcorner\left(\begin{matrix}
X \vspace{-4pt}\\ Z
\end{matrix}\right) = \urcorner(Z), \quad \llcorner(XY) =
\llcorner(Y),\quad \llcorner\left(\begin{matrix} X \vspace{-4pt}\\
Z \end{matrix}\right) = \llcorner(X).\end{align*}

\medbreak
\item[(d)] The double groupoid is vacant if and only if
the core groupoid is trivial.
\end{enumerate}

\newpage
\end{document}